\documentclass[11pt]{article}

\usepackage{amsxtra,amssymb,amsthm,amsmath,latexsym}
\usepackage{graphicx}
\usepackage{epsfig}
\usepackage{afterpage}

\newtheorem{lem}[]{Lemma}

 \newcommand{\lemref}[1]{Lemma~\ref{#1}}

\newcommand{\R}{{\mathbb R}}

\newcommand{\nb}{\nabla}

\newcommand{\dl}{{\delta}}
\newcommand{\bee}{\begin{equation*}}
\newcommand{\eee}{\end{equation*}}
\newcommand{\be}{\begin{equation}}
\newcommand{\ee}{\end{equation}}

\def\S_m{S}
\title{Electromagnetic Wave Scattering by  Small Impedance Particles
of an Arbitrary Shape}
\author{A.G. Ramm \\
\small Department of Mathematics\\[-0.8ex]
\small Kansas State University, Manhattan, KS 66506-2602, USA\\
\small \texttt{ramm@math.ksu.edu}}

\date{}
\begin{document}

% typeset front matter
\maketitle

\begin{abstract}
Scattering of electromagnetic (EM) waves by one and many small
($ka\ll 1$) impedance particles $D_m$ of an arbitrary shape,
embedded in a homogeneous medium, is studied. Analytic formula for
the field, scattered by one particle, is derived. The scattered
field is of the order $O(a^{2-\kappa})$, where $\kappa \in [0,1)$ is
a number. This field is much larger than in the Rayleigh-type
scattering.  An equation is derived for the effective EM field
scattered by many small impedance particles distributed in a bounded
domain. Novel physical effects in this domain are described and
discussed.

 \end{abstract}
%\pn{\\ {\em MSC:}\,\, 26D10;34G20; 37L05;44J05; 47J35;}




\section{Introduction}
In this paper we develop a theory of electromagnetic (EM) wave
scattering by one and many small  impedance particles (bodies)
$D_m$, $1\le m\le M=M(a)$, embedded in a homogeneous medium which is
described by the constant permittivity $\epsilon_0>0$ and
permeability $\mu_0>0$. The smallness of a particle means that
$ka\ll 1$, where $a=0.5 diam D_m$ is the characteristic dimension of
a particle, $k=\omega (\epsilon_0 \mu_0)^{1/2}$ is the wave number
in the medium exterior to the particles. Although scattering of EM
waves by small bodies has a long history, going back to Rayleigh
(1871), see \cite{LL}, \cite{R476}, the results of this paper are
new and useful in applications because light scattering by colloidal
particles in a solution, and light scattering by small dust
particles in the air are examples of the problems to which our
theory is applicable. The Mie theory deals with scattering by a
sphere, not necessarily small, and gives the solution to the
scattering problem in terms of the series in spherical harmonics. If
the sphere is small, $ka\ll 1$, then  the first term in the Mie
series yields the main part of the solution. Our theory is
applicable only to small particles. They can be of arbitrary shapes.
The solution to the scattering problem for one small particle
of an arbitrary shape  is
given  analytically. For many such particles the solution is reduced
to solving a linear algebraic system. This system is not obtained by
a discretization of some boundary integral equation, and it has a clear
physical meaning. Its limiting form as $a\to 0$ yields an
integro-differential equation for the limiting effective field in
the medium where the small particles are distributed. Wave
scattering problem in the case of one body can be studied
theoretically only in the limiting cases of a small body, $ka\ll 1$,
or a large body, $ka\gg 1$. In the latter case the geometrical optics
is applicable. This paper deals with the case $ka\ll 1$. Rayleigh
(1871) understood that the scattering by a small body is given
mainly by the dipole radiation. For a small body of arbitrary shape
this dipole radiation is determined by the polarization moment,
which is defined by the polarizability tensor. For homogeneous
bodies of arbitrary shapes analytical formulas, which allow one to
calculate this tensor with any desired accuracy, were derived in
\cite{R476}. These bodies were assumed dielectric or conducting in
\cite{R476}.  Under the Rayleigh assumption the scattered field is
proportional to $a^3$, to the volume of a small body. The physically
novel feature of our theory is a conclusion that for small impedance
particle with  impedance $\zeta=h(x)a^{-\kappa}$, where $h(x)$ is a
given function, Re$h(x)\ge 0$, and $\kappa\in [0, 1)$ is a given
number, the scattered field is proportional to $a^{2-\kappa}$, that
is, it is much larger than in the Rayleigh case, since
$a^3<<a^{2-\kappa}$ if $a<<1$ and $\kappa\ge 0$.

 In this paper wave scattering by small
impedance particles is studied. Besides high intrinsic interest in this
problem, the theory we develop allows one to get some physically
interesting conclusions about the changes of the material properties
of the medium in which many small particles are embedded. The results
of this paper can be used to develop a method for creating materials
with a desired refraction
coefficient by embedding many small impedance particles into a given
material. Such a theory has been developed by the author for scalar
wave scattering, for example, acoustic wave scattering, in a series
of papers \cite{R509}-\cite{R597}. The novel physical idea in this
paper is to reduce solving the scattering problem to finding some
constant pseudovector $Q$ (see formula \eqref{e25}), rather than a
pseudovector function $J$ (see formula \eqref{e8}) on the surface
of the scatterer. The quantity $Q$ is a pseudovector which somewhat
analogous to the total charge on the surface of the perfect conductor
with the shape of $D_m$,
while the function $J$ is analogous to the surface charge density.
We assume for simplicity that the impedance $\zeta$ (see formula
\eqref{e2}) is a constant given  in \eqref{e23}. A similar
assumption appeared in paper \cite{R536}, where
scalar wave scattering theory was developed. The results of the theory in
\cite{R536}  was a recipe for creating materials with a desired
refraction coefficient in acoustics (see \cite{R591}, \cite{R595},
\cite{R597}). The boundary impedance \eqref{e23}
grows to infinity as $a\to 0$, and one can pass to the limit
in the equation for the effective (self-consistent) field in the
medium, obtained by embedding many small impedance particles into a
given medium. Such a theory is briefly summarized in paper
\cite{R595}, where the equation for the limiting field in the medium
is derived. The aim of this paper is to develop a similar theory for
EM wave scattering by many small impedance particles embedded in a
given material.

For EM wave scattering by one small body $D$ of an arbitrary shape
with an impedance boundary condition an analytic formula for the
electromagnetic field in the region $r:=|x|\gg a$, is derived:
\be\label{e0} E(x)=E_0(x)+ [\nb \frac{e^{ikr}}{4\pi r},Q], \qquad
r\gg a, \qquad g(x,y):=\frac{e^{ikr}}{4\pi r},\quad r=|x-y|,\ee
where $E_0$ is the incident field, which satisfies Maxwell's
equations in the absence of the scatterer $D$, $[A,B]=A\times B$ is
the cross product of two vectors, $(Q,e_j)=Q\cdot e_j$ is the dot
product, $\{e_j\}_{j=1}^3$ is an orthonormal basis in $\R^3$,
\be\label{e0a} Q_j:=(Q,e_j)=-\frac{\zeta |S|}{i\omega
\mu_0}\Xi_{jp}(\nb\times E_0(O))_p, \quad \Xi:=(I+\alpha)\tau,\ee
over the repeated index $p$ summation is understood from $1$ to $3$,
$\zeta$ is the boundary impedance,
 $|S|$ is the surface area of the particle,  the matrix $\Xi_{jp}$ is
defined by the formula \be\label{e0b}
\Xi_{jp}:=(I+\alpha)\Big(\delta_{jp}-|S|^{-1}\int_S N_j(s)N_p(s)ds\Big):=
(I+\alpha)\tau_{jp},
\ee
where $N_j(s)$ is the $j-$th component of the unit normal $N(s)$ to
the surface $S$ at a point $s\in S$, pointing {\it out of} $D$,
$\tau:=I-b$, $b$ is a matrix, $b_{jp}:=|S|^{-1}\int_S N_j(s)N_p(s)ds$,
$k=\omega(\epsilon_0 \mu_0)^{1/2}$ is the wave number,  $O\in D$ is
the origin, $I$ is the identity matrix, and $\alpha$ is a matrix,
$$I+\alpha:=(I+\beta)^{-1},$$ where the matrix $\beta$ is defined in
\eqref{e35} (see below).

By $S^2$ the unit sphere in $\R^3$ will be denoted. The boundary $S$
of the small body $D$ is assumed smooth: it is sufficient to assume
that in local coordinates the equation of $S$ is given as $x_3=\phi
(x_1,x_2)$, where the function $\phi$ has first derivative
satisfying a H\"older condition.

Briefly speaking, there are three basic novel results in this
paper. The {first result} are formulas
\eqref{e0}-\eqref{e0b} and equations \eqref{e36a} and  \eqref{e37}.
 The {\it second result} is the reduction of the solution to
many-body scattering problem to solving a linear algebraic system
(LAS), see equations \eqref{e45} and  \eqref{e46}. The {\it third
result} is a derivation of the equation for
the limiting effective (self-consistent) field in the medium in
which many small impedance particles are embedded, see equation
\eqref{e49}.

 The scattering problem by one small body is
formulated and studied in Section 2, the reduction of the solution
to the many-body EM wave scattering problem to the solution of a LAS
is given in Section 3. Also in Section 3 a derivation of the
equation for the limiting effective field is given as $a\to 0$ and
$M=M(a)\to \infty$. In Section 4 the conclusions are formulated.

{\it In this paper we do not  solve the boundary integral equation
to which the scattering problem can be reduced in a standard
approach, but find asymptotically exact analytical expression for
the pseudovector $Q$ which defines the behavior of the scattered
field at distances $d>>a$.}

In fact, these distances $d$ can be very small if $a$ is
sufficiently small, and $d$ can be much less than the wavelength
$\lambda=\frac{2\pi}{k}$.

{\it Therefore our theory is valid in the physical situations in
which multiple scattering effects are dominant.}

\section{EM wave scattering by one small impedance particle}

Let  us use in this Section the following notations:  $D$ is a small
body, $D':=\R^3\setminus D,$ $k>0$ is a wave number, $ka\ll 1$,
$a=0.5 diam D$, $k=\frac {2\pi}{\lambda}$, $\lambda$ is the
wavelength of the incident EM wave, $k^2=\omega^2\epsilon_0\mu_0$,
where $\omega$ is the frequency and $\epsilon_0$, $\mu_0$ are
constant permittivity and permeability of the medium. Our arguments
remain valid if one assumes that the medium has a constant
conductivity $\sigma_0\ge 0$. In this case $\epsilon_0$ is replaced
by $\epsilon_0+i\frac {\sigma_0}{\omega}.$ Denote by $S$ the
boundary of $D$,  by $|S|$ its surface area, by $V$ the volume of
$D$, by $[E,H]=E\times H$ the cross product of two vectors, and by
$(E,H)=E\cdot H$ the dot product of two vectors, $N$ is the unit
normal to $S$ {\it pointing out of} $D$, $\zeta$ is the boundary
impedance of the particle.

Let $D$ be embedded in a homogeneous medium with constant parameters
$\epsilon_0$, $\mu_0$. Electromagnetic (EM) wave scattering problem
consists of finding vectors $E$ and $H$ satisfying Maxwell's
equations: \be\label{e1} \nb \times E=i\omega\mu_0 H,\quad \nb\times
H=-i\omega \epsilon_0 E\quad \text{in } D':=\R^3\setminus D, \ee the
impedance boundary condition:
\be\label{e2} [N,[E,N]]=\zeta
[N,H]\quad \text{ on } S \ee
and the radiation condition:
\be\label{e3} E=E_0+v_E,\quad H=H_0+v_H,
 \ee
where $E_0, H_0$ are the incident fields satisfying equations
\eqref{e1} in all of $\R^3$, $v_E=v$ and $v_H$ are the scattered
fields. In the literature, for example in \cite{LL}, the impedance
boundary condition is written sometimes as $E^t=\zeta[H^t,N]$, where
$N$ is the unit normal on $S$ pointed {\it into} $D$. Since our $N$
is pointed {\it out of} $D$, our impedance boundary condition \eqref{e2}
is the same as in \cite{LL}.

One often assumes that the incident wave is a plane wave, i.e.,
$E_0=\mathcal{E}e^{ik\alpha\cdot x}$, $\mathcal{E}$ is a constant
vector, $\alpha\in S^2$ is a unit vector, $S^2$ is the unit sphere
in $\R^3$, $\alpha\cdot \mathcal{E}=0$, $v_E $ and $v_H$ satisfy the
Sommerfeld radiation condition: $r(\frac{\partial v}{\partial
r}-ikv)=o(1)$ as $r:=|x|\to \infty$, and, consequently,
$[r^0,E]=H+o(r^{-1})$ and $[r^0,H]=-E+o(r^{-1})$ as $r:=|x|\to
\infty$, where $r^0:=x/r$.

{\it It is assumed in this paper that the impedance $\zeta$ is a
constant, Re $\zeta\ge 0$}.

This assumption guarantees the uniqueness of the solution to
Maxwell's equation satisfying the radiation condition. For
completeness a proof of the uniqueness result is given in Lemma 1.
 The tangential component of $E$ on $S$, $E^t$,
is defined as: \be\label{e5} E^t=E-N(E,N)=[N,[E,N]]. \ee This
definition differs from the one used often in the literature,
namely, from the definition $E^t=[N,E]$. Our definition \eqref{e5}
corresponds to the geometrical meaning of the tangential component
of $E$ and, therefore, should be used. The impedance boundary
condition is written in \cite{LL}  as $E^t=\zeta[H^t,N_i],$ where $\zeta$
is the  boundary impedance and $N_i$ is the unit normal to $S$ pointing
{\it into} $D$. In our paper $N$ is the unit normal pointing {\it out} of
$D$. Therefore, the impedance boundary condition in our paper is
written as in equation \eqref{e2}.  If one uses definition \eqref{e5},
then this condition reduces to \eqref{e2}, because $[[N,[H,N]],N]=[H,N].$ The
assumption Re$\zeta\ge 0$ is physically justified by the fact that
this assumption guarantees the uniqueness of the solution to the
boundary problem  \eqref{e1}-\eqref{e3}.

\begin{lem}\label{lem1}
Problem \eqref{e1}-\eqref{e3} has at most one solution.
\end{lem}
{\it Proof of \lemref{lem1}.} \noindent  In Lemma 1 it is not
assumed that $D$ is small. The proof is valid for an arbitrary
finite domain $D$. To prove the lemma one assumes that $E_0=H_0=0$
and has to prove that then $E=H=0$. Let the overbar stand for the
complex conjugate. The radiation condition implies
$[H,x^0]=E+o(|x|^{-1})$, where $x^0:=x/r$ is the unit normal on the
sphere  $S_R$ centered at the origin and of large radius $R$, and
one has $ I':=\int_{S_R} [E, \overline{H}]\cdot x^0 ds=\int_{S_R}
|E|^2ds+o(1)$ as $R\to \infty$.

From equations \eqref{e1} one derives:
$$\int_{D_{R}}(\overline{H}\cdot \nb\times E-E\cdot \nb
\times \overline{H})dx=\int_{D_R}(i\omega \mu_0 |H|^2-i\omega
\epsilon_0|E|^2)dx,$$ where $D_R:=D'\cap B_R$, and $R>0$ is so large that
$D\subset B_R:=\{x\ : \ |x|\leq R\}$. Recall that $\nb\cdot [E,
\overline{H}]=\overline{H}\cdot\nb\times E -E\cdot\nb\times
\overline{H}$. Applying the divergence theorem, using the radiation
condition on the sphere $S_R=\partial B_R$, and taking real part,
one gets
$$0=-\text{Re}\int_{S}[E,\overline{H}]\cdot N
ds+ \text{Re}\int_{S_R}[E,\overline{H}]\cdot x^0 ds:=I+I'.$$ The
radiation condition implies $I'\ge 0$ as $R\to \infty$.
The minus sign in front of the integral $I$ comes from the assumption that
$N$ on $S$ is directed out of $D$. The
impedance boundary condition and the assumption Re$ \zeta\ge 0$
implies $I\ge 0$. One has $I+I'=0$. One has
$$-\int_S[E,\overline{H}]\cdot Nds=\int_S E\cdot \overline{[N,H]}ds=
\int_SE\cdot \overline{[N,[E,N]]/\zeta},$$
and $E\cdot \overline{[N,[E,N]]/\zeta}=|E^t|^2\zeta|\zeta|^{-2}$.
Therefore, $I= \text{Re} \zeta
|\zeta|^{-2}\int_S |E_t|^2ds\ge 0$, and $I'\ge 0$, the
relation  $I+I'=0$ implies $I=0$. Therefore, if
$\text{Re}\zeta> 0$, then $E_t=0$ on $S$. Consequently, $E=H=0$ in
$D$. If $\text{Re}\zeta\ge 0$, then $\lim_{R\to \infty}\int_{S_R}|E|^2
ds=0$, and since $E$ satisfies the Sommerfeld radiation condition
it follows that $E=H=0$ in $D'$.\\
\lemref{lem1} is proved. \hfill $\Box$

 Let us note that problem \eqref{e1}-\eqref{e3} is equivalent to
the problem  \be\label{e6} \nb\times \nb\times E=k^2E\,\,\text{ in
}\,\, D',\quad H=\frac{\nb\times E}{i\omega \mu_0}, \ee
\be\label{e7} [N,[E,N]]=\frac{\zeta}{i\omega \mu_0}[N,\nb\times
E]\text{ on } S,\ee  together with the radiation condition
\eqref{e3}. Thus, we have reduced the scattering problem to finding
one vector $E(x)$. If $E(x)$ is found, then $H=\frac{\nb\times
E}{i\omega\mu_0},$ and the pair $E$ and $H$ solves Maxwell's
equations, satisfies the impedance boundary condition and the
radiation condition \eqref{e3}.

Let us look for $E$ of the form \be\label{e8} E=E_0+\nb \times
\int_{S}g(x,t)J(t)dt,\qquad g(x,y)=\frac{e^{ik|x-y|}}{4\pi|x-y|},
\ee where $E_0$ is the incident field, which satisfies Maxwell's
equations in the absence of the scatterer $D$, $t$ is a point on the
surface $S$, $t\in S$, $dt$ is an element of the area of $S$, and
$J(t)$ is an unknown pseudovector-function on $S$, which is
tangential to $S$, i.e., $N(t)\cdot J(t)=0$, where $N(t)$ is the
unit normal to $S$ at the point $t\in S$. That $J=J(t)$ is a
pseudovector follows from the fact that $E$ is a vector and $\nb
\times (gJ)$ is a vector only if $J$ is a pseudovector, because $g$ is a
scalar.

{\it It is assumed that $J$ is a smooth function on $S$, for
example, $J\in C^2(S)$.}

The right-hand side of \eqref{e8} solves equation \eqref{e6} in $D$ for
any continuous
$J(t)$, because $E_0$ solves \eqref{e6} and
\be\label{e9}\begin{split} &\nb\times\nb\times\nb\times
\int_{S}g(x,t)J(t)dt=grad div
\nb\times\int_{S}g(x,t)J(t)dt\\
&-\nb^2\nb\times \int_{S}g(x,t)J(t)dt =k^2\nb\times
\int_{S}g(x,t)J(t)dt,\quad x\in D'.
\end{split}\ee
Here we have used the known identity $div curl E=0,$ valid for any
smooth vector field $E$, and the known formula \be\label{eG} -\nb^2
g(x,y)=k^2g(x,y)+\dl(x-y). \ee The integral $\int_{S}g(x,t)J(t)dt$
satisfies the radiation condition. Thus, formula \eqref{e8} solves
problem \eqref{e6}, \eqref{e7}, \eqref{e3}, if $J(t)$ is chosen so
that boundary condition \eqref{e7} is satisfied.

Let $O\in D$ be a point inside $D$. The following  known formula
(see, for example, \cite{M}) is useful: \be\label{e12} [N,\nb\times
\int_{S}g(x,t)J(t)dt]_{\mp}=\int_{S}
[N_s,[\nb_xg(x,t)|_{x=s},J(t)]]dt\pm \frac{J(s)}{2}, \ee where the
$\pm$ signs denote the limiting values of the left-hand side of
\eqref{e12} as $x\to s$ from $D$, respectively, from  $D'$. To
derive an integral equation for $J=J(t)$, substitute $E(x)$ from
\eqref{e8} into impedance boundary condition \eqref{e7},
 and get
 \be\label{e13}
0=f+[N,[\nb\times\int_Sg(s,t)J(t)dt, N]] -\frac{\zeta}{i\omega
\mu_0}[ N, \nb\times \nb \times \int_Sg(s,t)J(t)dt],\ee where
 \be\label{e14}
f:=[N,[E_0(s),N]]-\frac{\zeta}{i\omega \mu_0}[N,\nb\times E_0].
\ee

We  assume that \be\label{e23} \zeta=\frac{h}{a^\kappa},  \ee where
Re $h\ge 0$. and $\kappa \in [0,1)$ is a constant.

Let us write \eqref{e8} as \be\label{e24} E(x)=E_0(x)+[\nb_x
g(x,O),Q]+ \nb\times\int_{S}(g(x,t)-g(x,O))J(t)dt, \ee where
\be\label{e25} Q:=\int_{S}J(t)dt. \ee The central physical idea of
the theory, developed in this paper, is: the third term in the
right-hand side of \eqref{e24} is negligible compared with the
second term if $ka\ll 1$. Consequently, the scattering problem is
solved if $Q$ is found. A traditional approach requires finding an
unknown function $J(t)$, which is usually found numerically by the
boundary integral equations (BIE) method. The reason for the third
term in the right-hand side of \eqref{e24} to be negligible compared
with the second one, is explained by the estimates, given below. In
these estimates the smallness of the body is used essentially: even
if one is in the far zone, i.e., $\frac {a}{ d}\ll 1$, one cannot
conclude that estimate \eqref{e28} (see below) holds unless one
assumes that $ka\ll 1$. Thus, {\it the second sum in \eqref{e24}
cannot be neglected in the far zone if the condition $ka \ll 1$ does
not hold.}

We prove below that $Q=O(a^{2-\kappa})$. To prove that the third
term in the right-hand side of \eqref{e24} is negligible compared
with the second one, let us establish  several estimates valid if
$a\to 0$ and $d:=|x-O|\gg a$. Under these assumptions one has
\be\label{e26} j_1:=|[\nb_x g(x,O),Q]|\leq
O\left(\max\left\{\frac{1}{d^2},\frac{k}{d}\right\}\right)O(a^{2-\kappa}),
\ee \be\label{e27} j_2:=|\nb\times\int_{S}(g(x,t)-g(x,O))\sigma
(t)dt|\leq a
O\left(\max\left\{\frac{1}{d^3},\frac{k^2}{d}\right\}\right)O(a^{2-\kappa}),
\ee and \be\label{e28} \left| \frac{j_2}{j_1}\right|=O\left( \max
\left\{ \frac{a}{d},ka\right\}\right)\to 0,\qquad \frac a d=o(1),
\qquad a\to 0.  \ee {\it These estimates show that one may neglect
the third term in \eqref{e24}, and write} \be\label{e29}
E(x)=E_0(x)+[\nb_xg(x,O),Q]. \ee The error of this formula tends to
zero as $a\to 0$ under our assumptions.

Note that the inequality $|x|\gg ka^2$ is satisfied for $d=O(a)$ if
$ka\ll 1$. Thus, formula \eqref{e29} is applicable in a wide region.

{\it Let us estimate $Q$
asymptotically, as $a\to 0$.}

Take cross product of \eqref{e13} with $N$ and integrate the
resulting equation over $S$ to get
$$ I_0+I_1+I_2=0,$$
 where $I_0$ is
defined in \eqref{e30} (see below), $I_1$ is defined in
\eqref{e32a}, and
$$I_2:=-\int_S\frac{\zeta}{i\omega \mu_0}[N,[ N, \nb\times \nb
\times \int_Sg(s,t)J(t)dt]] ds.$$ Let us estimate the order of $I_0$ as
$a\to 0$. One has
 \be\label{e30} I_0=\int_S\Big([N,E_0]-\frac
{\zeta}{i\omega \mu_0}[N,[N,\nb \times E_0]]\Big)ds.\ee In what
follows we keep only the main terms as $a\to 0$, and denote by the
sign $\simeq$ the terms equivalent up to the terms of higher order
of smallness as $a\to 0$. One has
 $$I_{00}:=\int_S [N,E_0]ds=\int_D\nb \times E_0 dx\simeq \nb \times
E_0(O)V,$$ where $dx$ is the element of the volume,  $O\in D$ is a
point chosen as the origin, and $V$ is the volume of $D$. Thus,
$I_{00}=O(a^3)$. Denoting by $|S|$ the surface area of $S$, one
obtains \be\label{e31}I_{01}:=-\int_S \frac {\zeta}{i\omega
\mu_0}[N,[N,\nb \times E_0]]ds=\frac {\zeta |S|}{i\omega \mu_0}\tau
\nb \times E_0,\ee where \be\label{e32} \tau:=I-b, \quad
b_{mj}:=\frac 1 {|S|}\int_S N_m(t)N_j(t)dt,\ee $b=(b_{mj})$ is a
matrix which depends on the shape of $S$, and $I:=\delta_{mj}$ is
the unit matrix. Since $\zeta=O(a^{-\kappa})$ one concludes that
$I_{01}=O(a^{2-\kappa})$, so $|I_{00}|\ll |I_{01}|$, because
$I_{00}=O(a^3)$ as $a\to 0$. Thus,
$$I_0\simeq I_{01}=O(a^{2-\kappa}), \quad a\to 0.$$

 Let us consider $I_1$: \be\label{e32a} I_1=\int_S\Big(\int_S[N,[\nb
 g(s,t),J(t)]]dt +\frac {J(s)}{2}\Big)ds:=I_{11}+\frac Q 2.\ee
One has
$$I_{11}=\int_Sds \int_Sdt\Big(\nb g(s,t)N(s)\cdot J(t)- J(t) \frac
{\partial g(s,t)}{\partial N(s)}\Big).$$ It is well known that
$\int_S\frac{\partial g_0(s,t)}{\partial N(s)}ds=-\frac 1 2,$ where
$g_0:=\frac 1 {4\pi|s-t|}$. Since
$$g(s,t)-g_0(s,t)=\frac
{ik}{4\pi}+O(|s-t|), \quad  |s-t|\to 0,$$ and $|s-t|=O(a)$, one
concludes that if $a\to 0$ then \be\label{e33} \int_S \frac{\partial
g(s,t)}{\partial N(s)}ds\simeq -\frac 1 2,\ee and
\be\label{e34}\int_Sdt \int_Sds\nb g(s,t)N(s)\cdot J(t):=\beta Q,\ee
where the matrix $\beta$ is defined  by the formula: \be\label{e35}
\beta:=(\beta_{mj}) := \int_S \frac {\partial g(s,t)}{\partial
s_m}N_j(s)ds.\ee Therefore, \be\label{e36} I_1\simeq (I+\beta)Q.\ee
Matrices $\beta$ and $b$ for spheres are calculated at the end of
Section 1.

Let us show that $I_2=O(a^{4-\kappa})$ and therefore $I_2$
is negligible compared with $I_0$ as $a\to 0$. If this is done, then
equation for $Q$ is \be\label{e36a} Q=- \frac {\zeta |S|}{i\omega
\mu_0} (I+\alpha)\tau \nb \times E_0,\ee where
$$I+\alpha:=(I+\beta)^{-1},$$ and the matrix $I+\beta$ is invertible.
From \eqref{e36a} it follows that
 \be\label{e37} E(x)=E_0(x)-\frac{
\zeta |S|}{i\omega \mu_0}
 [\nb_x  g(x,O), (I+\alpha)\tau \nb\times E_0(O)]. \ee
{\it This equation is our first main result which gives an analytic
formula for the solution of the EM wave scattering problem by a small body
of an arbitrary shape, on the boundary of which an impedance boundary
condition holds.}

In the far zone $r:=|x|\to \infty$ one has $\nb_x
g(x,O)=ikg(x,O)x^0+O(r^{-2})$, where $x^0:=x/r$ is a unit vector in
the direction of $x$. Consequently, for $r\to \infty$ one can
rewrite formula \eqref{e37} as
\be\label{e37a} E(x)=E_0(x)-\zeta
|S|\Big(\frac{\epsilon_0}{\mu_0}\Big)^{1/2}
 \frac{e^{ikr}}{r}[x^0, (I+\alpha)\tau \nb\times E_0(O)]. \ee
This field is orthogonal to the radius-vector $x$ in the far zone as
it should be.

Let us show that the term $I_2$ is negligible as $a\to 0$. Remember that
$curl curl =grad div -\nb^2$ and
$$-\nb^2 \int_S
g(x,t)J(t)dt=k^2\int_S g(x,t)J(t)dt.$$
 Consequently,
$$-i\omega\mu_0I_2\simeq \zeta \int_S ds[N,[N, grad div \int_S
g(x,t)J(t)dt]]|_{x\to s}:=I_{21}.$$ Since the function $J(t)$ is
assumed smooth, one has
$$ div \int_S g(x,t)J(t)dt=\int_S
g(x,t)div J(t)dt,\quad div J=\frac {\partial J_m}{\partial t_m},$$
summation is understood here and below over the repeated indices,
and $div J$, where $J$ is a tangential to $S$ field, is the surface
divergence. Furthermore,
$$grad \int_S g(x,t)\frac {\partial J_m}{\partial t_m}dt=
e_p \int_S g(x,t)\frac {\partial^2 J_m}{\partial t_p \partial
t_m}dt, $$
where the relation $\frac {\partial g(x,t)}{\partial x_p}=
-\frac {\partial g(x,t)}{\partial t_p}$ was used, and an integration
by parts with respect to $t_p$ has been done over the closed surface
$S$. Therefore
$$I_2\le c|\zeta||\int_Sds
\int_S|g(s,t)|=O(a^{4-\kappa})\ll I_0.$$ The constant $c>0$ here is a bound
on the second derivatives of $J$ on $S$.

{\it Example of calculation of matrices $\beta$ and $b$}.

Let us calculate $\beta$ and $b$ for a sphere of radius $a$ centered
at the origin. It is quite easy to calculate $b$:
$$b_{jm}=\frac 1{|S|}\int_SN_jN_mdt=\frac 1 3 \delta_{jm}.$$
Note that $|S|=4\pi a^2$, $dt=a^2\sin \theta d\theta d\phi$, and
$N_m$ is proportional to the spherical function $Y_{1,m}$, so the
above formula for $j\neq m$ follows from the orthogonality
properties of the spherical functions, and for $m=j$ this formula is
a consequence of the normalization $|N|=1$.

It is less simple to calculate matrix $\beta$. One has
$$\beta_{jm}= \frac
{1}{4\pi}\int_0^{2\pi}d\phi \int_0^\pi d\theta \frac{\sin \theta
(t_j-s_j)s_m}{|s-t|^3}.$$ By symmetry one can choose $t=(0,0,1)$ and
let the $z-$axis pass through $t$ from the origin. Then
$|t-s|^2=2-2\cos\theta$, $1-\cos\theta=2\sin^2(\theta/2)$, and
$$\beta_{jm}= \frac
1{4\pi}\int_0^{2\pi}d\phi \int_0^\pi d\theta \frac{\sin \theta
(s_m\delta_{3j}-s_js_m)}{8 \sin^3(\theta/2)}.$$
This can be written
as $$\beta_{jm}= \frac 1{16\pi}\int_0^{2\pi}d\phi \int_0^\pi d\theta
\frac{\cos( \theta/2) (s_m\delta_{3j}-s_js_m)}{ \sin^2(\theta/2)}.$$
If $j\neq m$, then $\beta_{jm}=0$, as one can check. If $j=m\le 2$,
then $\delta_{j3}=0$ and $\beta_{jj}=\frac 1 6$, as one can check.
Finally,
$$\beta_{33}=\frac 1{16\pi}\int_0^{2\pi}d\phi \int_0^\pi
d\theta \cos( \theta/2) \frac {\cos \theta -\cos^2\theta)}{
\sin^2(\theta/2)}=\frac 1 6.$$  Thus,
$$\beta_{jm}=\frac 1 6 \delta_{jm}.$$

Let us formulate the result of this Section as a theorem.

{\bf Theorem 1.} {\it If $ka\ll 1$, then the solution to the
scattering problem \eqref{e1}-\eqref{e3} is given by formula
\eqref{e37}.}

In the next Section the EM wave scattering problem is studied in the
case of many small bodies (particles) whose physical properties are
described by the boundary impedance conditions.

\section{ Many-body EM wave scattering problem}
Consider now EM wave scattering by many small bodies (particles)
$D_m$, $1\le m \le M=M(a)$, distributed in an arbitrary bounded
domain $\Omega$.  Here $a$ is the characteristic dimension of a
particle, $M(a)$ is specified below, and $x_m\in D_m$ is a point.
Assume for simplicity of formulation of the problem that the
particles are of the same shape and orientation, and their physical
properties are described by the boundary impedance \eqref{e23},
where $h=h(x_m)$, and $h(x)$, Re$h(x)\ge 0$, is a continuous
function in an arbitrary large but finite domain $\Omega$ where the
particles are distributed. For simplicity let us assume that the
particles are distributed in a homogeneous medium with parameters
$\epsilon_0, \mu_0$.

The distribution of small particles is given by the following law:
\be\label{e38}
\mathcal{N}(\Delta)=a^{-(2-\kappa)}\int_{\Delta}N(x)dx \big(
1+o(1)\big), \ee where $\Delta$ is an arbitrary open subset of
$\Omega$, $N(x)\ge 0$ is a given function which one can choose as
one wishes, $\mathcal{N}(\Delta)$ is the number of particles in
$\Delta$, and $\kappa \in (0,1)$ is a parameter from formula
\eqref{e23}. The distance $d=d(a)\gg a$ between the neighboring
particles is assumed to satisfy the relation $\lim_{a\to
0}ad^{-1}(a)=0$. The order of magnitude of $d$ as $a\to 0$ can be
found by a simple argument: if one assumes that the particles are
packed so that between each neighboring pair of particles the
distance is $d$, then in a unit cube $\Delta_1$ the number of
particles is $O(d^{-3})$, and this amount should be equal to
$\mathcal{N}(\Delta_1)=O(a^{-(2-\kappa)})$. Therefore
$$d=O(a^{(2-\kappa)/3}).$$ Denote $D:=\bigcup_{m=1}^M D_m$,
$D':=R^3\setminus D$.

The scattering problem consists of finding the solution to
 \be\label{e39}\nb \times \nb \times E=k^2E\,\,\,  in \,\,\, D',\ee
\be\label{e40} [N,[E,N]]-\zeta_m [N, \nb \times E] \,\,\, on \,\,\,
S_m, \,\,\, 1\le m \le M,\ee \be\label{e41} E=E_0+v_E,\ee where
$E_0$ is the incident field, which satisfies equation \eqref{e39} in
$\R^3$, and  $v_E$ satisfies the radiation condition. Existence and
uniqueness of the solution to the scattering problem
\eqref{e39}-\eqref{e41} is known. Let us look for the solution of
the form \be\label{e42} E(x)=E_0(x)+ \sum_{m=1}^M \nb \times
\int_{S_m}g(x,t)J_m(t)dt.\ee {\it Define the effective
(self-consistent) field, acting on the $j-$th particle, by the
formula} \be\label{e43}
 E_e(x)=E_0(x)+ \sum_{m=1, m\neq j}^M \nb \times
\int_{S_m}g(x,t)J_m(t)dt.\ee Because of our assumption $d\gg a$, any
single particle $D_j$ can be considered as being placed in the
incident field $E_e(x)$. Therefore \be\label{e44}  E_e(x)=E_0(x)-
\sum_{m=1, m\neq j}^M \frac{ \zeta_m |S|}{i\omega \mu_0}
 [\nb_x  g(x,x_m), (I+\alpha)\tau \nb\times E_e(x_m)].
\ee Note that $\zeta |S|=h(x_m)c_Sa^{2-\kappa}$, where the constant
$c_S$ depends on the shape of $S$ and is defined by the formula
$$|S|=c_Sa^2.$$
 Thus, formula \eqref{e44} implies
 \be\label{e45}  E_e(x)=E_0(x)- \frac{c_S}{i\omega \mu_0}
a^{2-\kappa}\sum_{m=1,
m\neq j}^M h(x_m) [\nb_x g(x,x_m), (I+\alpha)\tau \nb\times
E_e(x_m)]. \ee Denote $A_m:=\nb\times E_e(x_m)$. Applying operator
$\nb \times$ to \eqref{e45} and then setting $x=x_j$, one obtains a
linear algebraic system (LAS) for the unknown $A_m$: \be\label{e46}
A_j=A_{0j}-\frac{c_S}{i\omega \mu_0}  a^{2-\kappa}\sum_{m=1, m\neq
j}^M h(x_m)  \nb_x \times [\nb_x g(x,x_m)|, (I+\alpha)\tau
A_m]|_{x=x_j},\ee
where $1\le j\le M$. If $A_m$ are found, then formula
\eqref{e44} yields the solution to many-body EM wave scattering
problem. The magnetic field is calculated by the formula $H(x)=\frac
{\nb \times E(x)}{i\omega \mu_0}$, if $E(x)$ is found.

{\it Let us now pass to the limit in equation \eqref{e45} as $a\to
0$.}

Consider a partition of $\Omega$ into a union of cubes $\Delta_p$,
$1\le p\le P$. These cubes have no common interior points, and their
side $\ell\gg d$, $\ell=\ell(a)\to 0$ as $a\to 0$. Choose a point
$x_p\in \Delta_p$, for example, let $x_p$ be the center of
$\Delta_p$. Rewrite \eqref{e45} as \be\label{e47}\begin{split}
 &E_e(x_q)=E_0(x_q)- \frac{c_S}{i\omega \mu_0} \sum_{p=1,
p\neq q}^P h(x_p) [\nb_x g(x,x_p)|_{x=x_q}, (I+\alpha)\tau \nb\times
E_e(x_p)]\\
&\cdot a^{2-\kappa} \sum_{x_m\in \Delta_p} 1.
\end{split}\ee
It follows from equation \eqref{e38} that
$$a^{2-\kappa}\sum_{x_m\in \Delta_p} 1\simeq N(x_p)|\Delta_p|,\quad a\to 0,$$
where $|\Delta_p|$ denotes the volume of the cube $\Delta_p$. Thus,
\eqref{e47} can be written as: \be\label{e48}\begin{split} &
E_e(x_q)=E_0(x_q)- \\
&\frac{c_S}{i\omega \mu_0}  \sum_{p=1, p\neq q}^P h(x_p)N(x_p)
[\nb_x g(x,x_p)|_{x=x_q}, (I+\alpha)\tau \nb\times
E_e(x_p)]|\Delta_p|.
\end{split}\ee This equation is a discretized version of
the integral equation:
 \be\label{e49} E(x)=E_0(x)- \frac{c_S}{i\omega \mu_0}
 \nb \times \int_{\Omega}  g(x,y) (I+\alpha)\tau
\nb\times E(y) N(y)h(y)dy.\ee Equation \eqref{e49} is the equation
for the limiting, as $a\to 0$, effective field in the medium in
which $M(a)$ small impedance particles are embedded according to
 the distribution law \eqref{e38}.

Let us summarize our result.

{\bf Theorem 2.} {\it If small particles are distributed according
to \eqref{e38} and their boundary impedances are defined in
\eqref{e23}, then the solution to a many-body EM wave scattering
problem \eqref{e39}-\eqref{e41} is given by formula \eqref{e44},
where the quantities $\nb\times E(x_m):= A_m$ are found from the
linear algebraic system \eqref{e46}, and the limiting, as $a\to 0$,
electric field $E(x)$ in $\Omega$ solves integral equation
\eqref{e49}.}

 Let us discuss the {\it novel physical consequences  of our theory.}
Applying the operator $curl curl $ to equation \eqref{e49} and
taking into account that $curl curl=grad div- \Delta$, $div
curl\equiv 0$, and $-\Delta g(x,y)=k^2g(x,y)+\delta(x-y)$ one
obtains \be\label{e50} \nb \times \nb \times E=k^2E -\nb \times
 \Big(\frac{c_S}{i\omega \mu_0}(I+\alpha)\tau \nb\times E(x)
N(x)h(x)\Big).\ee The novel term, which is  due to the limiting
distribution of the small particles, is $$T:=-\nb \times
\Big(\frac{c_S}{i\omega \mu_0}(I+\alpha)\tau \nb\times E(x)
N(x)h(x)\Big).$$ Let us make a simple assumption that
$N(x)h(x)=const$. Physically this means that the distribution of the
small particles is uniform in the space and the boundary impedances
of the particles do not vary in the space. Denote
$\frac{c_S}{i\omega \mu_0}Nh:=c_1$. Assume also that $\alpha$ and
$\tau$ are proportional to the identity matrix. This, for example,
happens if all $S_m$ are spheres. Then $T=-c_2\nb\times \nb \times
E(x)$, where $c_2=const$.

{\it In this case the novel term  $T$ can be interpreted physically
in a simple way: it results to a change in the refraction
coefficient of the medium in $\Omega$.}

Indeed, by taking the term $T$ to the left side of \eqref{e50} and
dividing both sides of this equation by $1+c_2$, one sees that the
coefficient $k^2$ is replaced by $k_1^2:=\frac{k^2}{1+c_2}$. Since
$c_2$ is a complex number, the new medium is absorptive. The
assumption Re$h\le 0$ implies that Im$k_1^2\ge 0$. Since $h(x)$ and
$N(x)$  can be chosen by the experimentalist  as he wishes, their
dependence on $x$ can be chosen as he wishes. Therefore,
distributing many small particles in $\Omega$ one can change the
refraction coefficient of the medium in a desirable direction.

\section{Conclusions}
The main results of this paper are:

1. Equations \eqref{e37} and \eqref{e36a}
for the solution of the EM wave scattering
problem \eqref{e1}-\eqref{e3} by one small body of an arbitrary
shape.

2. A numerical method, based on  linear algebraic system (LAS)
\eqref{e46}, see formulas \eqref{e45} and \eqref{e46}, for solving
many-body EM wave scattering problem in the
case of small bodies of an arbitrary shape.

3. Equation \eqref{e49} for the electric field in the limiting
medium obtained by embedding $M(a)$ small impedance particles of an
arbitrary shape, distributed according to \eqref{e38},  as $a\to 0$.

%\newpage

\end{document}